\documentclass[12pt]{article}
\usepackage{amssymb}
\newcommand{\RR}{\mathbb{R}}
\newcommand{\KK}{\mathbb{K}}

\newcommand{\PP}{\mathbb{P}}
\newcommand{\QQ}{\mathbb{Q}}

\newcommand{\CC}{\mathbb{C}}
\newcommand{\DD}{\mathbb{D}}
\newcommand{\NN}{\mathbb{N}}

\newcommand{\sphere}{\mathbb{S}}
\def\build#1_#2^#3{\mathrel{
\mathop{\kern 0pt#1}\limits_{#2}^{#3}}}
\def\cq{$\hfill \square$}
\def \un{\underline}
\def\t{{\cal T}}
\def\W{{\cal W}}
\def\w{{\rm w}}

\def\eg{{\bf e}}
\def\ov{\overline}
\def\wh{\widehat}
\def\wt{\widetilde}
\def\la{\longrightarrow}

\def\noi{\noindent}
\def\dem{\vskip 3mm \noindent{\bf Proof.}\hskip10pt}
\def\rem{\noindent{\bf Remark.} }

\newtheorem{theorem}{Theorem}[section]
\newtheorem{lemma}[theorem]{Lemma}
\newtheorem{proposition}[theorem]{Proposition}
\newtheorem{corollary}[theorem]{Corollary}

\title{Scaling limits of bipartite planar maps \\ are
homeomorphic to the $2$-sphere}
\author{Jean-Fran\c cois Le Gall \and Fr\'ed\'eric Paulin}
\date{\small\today}

\begin{document}
\maketitle

\begin{abstract}
  We prove that scaling limits of random planar maps which are
  uniformly distributed over the set of all rooted $2k$-angulations
  are a.s.~homeomorphic to the two-dimensional sphere. Our methods
  rely on the study of certain random geodesic laminations of the
  disk.
\end{abstract}

\section{Introduction}
\label{sec:intro}

This paper continues the study of scaling limits of large random
planar maps in the sense of the Hausdorff-Gromov topology. In the
particular case of uniformly distributed $2k$-angulations, scaling
limits were shown in \cite{L} to be homeomorphic to a (random) compact
metric space which may be naturally defined as a quotient of the
Continuum Random Tree (CRT), which was introduced by Aldous in
\cite{Al1,Al3}.  The main goal of the present paper is to prove that
this limiting metric space is almost surely homeomorphic to the
$2$-sphere $\sphere^2$.

Let us first recall some basic definitions. More details can be found
in \cite{L}. A {\it planar map} is a topological embedding (without
edge crossing) of a finite connected graph in the sphere $\sphere^2$. Its
{\it faces} are the connected components of the complement of its
image in $\sphere^2$. Let $k\geq 2$ be a fixed integer.  A $2k$-angulation
is a planar map such that each face is adjacent to $2k$ edges (one
should in fact count edge sides, so that if an edge lies entirely
inside a face, it should be counted twice). A planar map is called
{\it rooted} if it has a distinguished oriented edge, which is called
the {\it root edge}. Two rooted planar maps are said to be {\it
equivalent} if the second one is the image of the first one under an
orientation-preserving homeomorphism of the sphere, which also
preserves the root edge. We systematically identify equivalent rooted
planar maps. Thanks to this identification, the set of all rooted
$2k$-angulations with a given number of faces is finite.

For every integer $n\geq 2$, let ${\cal M}^k_n$ be the set of all
rooted $2k$-angulations with $n$ faces, and let $M_n$ be a random
planar map that is uniformly distributed over ${\cal M}^k_n$. Denote
by ${\bf m}_n$ the set of vertices of $M_n$, and write $d_n$ for the
graph distance on ${\bf m}_n$. We view $({\bf m}_n,d_n)$ as a
random compact metric space, and study its convergence in
distribution as $n\la\infty$, after a suitable rescaling.

We denote by $\KK$ the set of all isometry classes of compact metric
spaces, and equip $\KK$ with the Hausdorff-Gromov distance $d_{GH}$
(see \cite{Gro}, \cite{P} or \cite{BBI}). Then $(\KK,d_{GH})$ is a
Polish space, which makes it appropriate to study the convergence in
distribution of $\KK$-valued random variables.

We can now state our main result.

\begin{theorem}
\label{theo:intromain} %
The sequence of the laws of the metric spaces $({\bf m}_n,
n^{-1/4}d_n)$ is tight (i.e. relatively compact) in the space of all
probability measures on $\KK$. If $({\bf m}_\infty,d_\infty)$ is the
weak limit of a subsequence of $({\bf m}_n,n^{-1/4}d_n)$, then the
metric space $({\bf m}_\infty,d_\infty)$ is almost surely homeomorphic
to the sphere $\sphere^2$.
\end{theorem}

\rem %
It is natural to conjecture that the sequence $({\bf m}_n, n^{-1/4}
d_n)$ does converge in distribution, or equivalently that the law of
any weak limit $({\bf m}_\infty,d_\infty)$ is uniquely determined, and
that this law is independent of $k$ up to multiplicative constants.
This is still an open problem, even though detailed information on
$({\bf m}_\infty,d_\infty)$ is already available.  In particular, it
is known that the Hausdorff dimension of $({\bf m}_\infty, d_\infty)$
is almost surely equal to $4$ (see \cite[Theorem 6.1]{L}).

\medskip %
The first assertion of Theorem \ref{theo:intromain} is already stated
in Proposition 3.2 of \cite{L}. The new part of the theorem is the
second assertion, which is proved in Section \ref{sec:proofmain}
below. We rely on the main theorem of \cite{L}, which asserts in
particular that any weak limit $({\bf m}_\infty,d_\infty)$ is almost
surely homeomorphic to a quotient of the CRT corresponding to a
certain pseudo-metric $D^*$ (see Section \ref{sec:proofmain} for
details).  As a preparation for the proof of our main result, Section
\ref{sec:treegeodlam} investigates, in a deterministic setting,
quotient spaces of compact $\RR$-trees coded by continuous functions
on the circle, and their relations with geodesic laminations of the
disk.  As a matter of fact, a key idea is to observe that the CRT,
which is the random $\RR$-tree coded by a normalized Brownian
excursion (in the sense of Theorem 2.1 of \cite{DuLG}), can also be
interpreted as the quotient space induced by a certain random geodesic
lamination of the hyperbolic disk. This observation is related to the
work of Aldous \cite{Al4,Altri} about random triangulations of the
circle: The random geodesic lamination that we consider corresponds to
the random triangulation in Section 5 of \cite{Al4} (or Section 2.3 in
\cite{Altri}), provided we replace the Poincar\'e disk model of
Lobatchevsky's hyperbolic plane with the Klein disk model.

Any random metric space that arises as a weak limit of rescaled planar
maps is then homeomorphic to a topological space that can be obtained
by taking one more quotient with respect to a second random geodesic
lamination, which is not independent of the first one.  To handle this
setting, we introduce on the sphere $\sphere^2$ the equivalence relation
for which two distinct points of the upper hemisphere, resp.~of the
lower hemisphere, are equivalent if they belong to the same geodesic
line of the first random lamination, resp.~of the second one, or to
the closure of an ideal hyperbolic triangle which is a connected
component of the complement of the same lamination. To get the second
assertion of Theorem \ref{theo:intromain}, we then use a theorem of
Moore \cite{Moo} giving sufficient conditions for a quotient space of
the sphere $\sphere^2$ to be homeomorphic to the sphere.

\medskip %
Theorem \ref{theo:intromain} yields information about the large scale
geometry of random planar maps. Let us state a typical result in this
direction.  Recall that a path of length $p$ in a planar map is a
sequence $x_0,e_1,x_1,e_2,\ldots,x_{p-1},e_{p},x_p$, where
$x_0,x_1,\ldots,x_p$ are vertices, $e_1,\ldots,e_p$ are edges and the
endpoints of $e_i$ are the points $x_{i-1}$ and $x_i$, for every
$i\in\{1,\ldots,p\}$.  The path is called a cycle if $x_0=x_p$. We say
that it is an injective cycle if in addition $x_1,\ldots,x_p$ are
distinct (when $p=2$, we also require that $e_1\ne e_2$). If $C$ is an
injective cycle, then the union of its edges ${\cal R}(C)$ separates
the sphere in two connected components, by Jordan's theorem.

\begin{corollary}
\label{coro:introtype} %
Let $\delta>0$ and let $\theta:\NN\la \RR_+$ be a function such that
$\theta(n)=o(n^{1/4})$ as $n\to\infty$. Then, with a probability
tending to $1$ as $n\to \infty$, there exists no injective cycle $C$
of the map $M_n$ with length $\ell(C)\leq \theta(n)$ such that the set
of vertices that lie in either connected component of $\sphere^2\setminus
{\cal R}(C)$ has diameter at least $\delta n^{1/4}$.
\end{corollary}

Notice that the diameter of the map $M_n$ is of order $n^{1/4}$ by
Theorem \ref{theo:intromain} (see also Theorem 3 in \cite{MaMi} or
Theorem 2.5 in \cite{We}). So Corollary \ref{coro:introtype} says that
with a probability close to one when $n\to\infty$, we cannot find
small ``bottlenecks'' in the map $M_n$ such that both sides of the
bottleneck have a diameter which is also of order $n^{1/4}$.

We refer to the introduction of \cite{L} for a detailed discussion of
the recent work about asymptotics for random planar maps. The idea of
studying the scaling limit of random quadrangulations appeared in
Chassaing and Schaeffer \cite{CS}.  This paper made an extensive use
of bijections between quadrangulations and trees, which have been
extended to very general planar maps by Bouttier, Di Francesco and
Guitter \cite{BDG}.  Marckert and Mokkadem \cite{MaMo} conjectured
that the scaling limit of random quadrangulations should be given by
the so-called Brownian map, which is essentially the same object as
the quotient of the CRT that was mentioned above (see also \cite{MaMi}
for related work on more general planar maps). Planar maps play an
important role in theoretical physics. See the pioneering paper
\cite{BIPZ} for the relation between enumeration problems for maps and
the evaluation of matrix integrals. Bouttier's thesis \cite{Bo} gives
an overview of the connections between planar maps and statistical
physics.

As a final remark, it is very likely that Theorem \ref{theo:intromain}
can be extended to more general random planar maps, in particular to
the Boltzmann distributions on bipartite planar maps which are
discussed in \cite{MaMi}.  The recent work of Miermont \cite{Mir} also
suggests that similar results should hold for random triangulations.

The paper is organized as follows. Section 2 introduces the $\RR$-tree
$\t_g$ coded by a continuous function $g$ on the circle, and
associates with this tree a geodesic lamination $L_g$ of the disk.
Moore's theorem is used in the proof of Proposition
\ref{prop:quotspher} to verify that certain quotients of $\t_g$ are
homeomorphic to the sphere $\sphere^2$. In addition, Section 2 gives a few
properties of the lamination $L_g$, and in particular computes its
Hausdorff dimension under suitable assumptions on the function $g$
(Proposition \ref{Hausdorff-lami}).  In the particular case when $g$
is the normalized Brownian excursion, one recovers the value $3/2$
which was given in \cite{Al4} (see Proposition \ref{Hausdorff-random}
below). Section 3 contains the proof of our main results. The key step
is to verify that any weak limit in Theorem \ref{theo:intromain} can
be written in the form of a quotient space which satisfies the
assumptions needed to apply Proposition \ref{prop:quotspher}.  The
verification of these assumptions requires two technical lemmas, whose
proofs are postponed to Section 4. The path-valued random process
called the Brownian snake plays an important role in these proofs.

\medskip
{\small \noi{\bf Acknowledgement.} We thank Andrei Okounkov for a
useful remark to the first author that motivated Corollary
\ref{coro:introtype}.}

\section{Trees and geodesic laminations}
\label{sec:treegeodlam}

In this section, we deal with various quotient spaces. Let $E$ be a
topological space, and let $\sim$ be an equivalence relation on $E$.
Unless otherwise stated, the quotient space $E/\!\sim$ will always be
equipped with the quotient topology, which is the finest topology on
$E/\!\sim$ such that the canonical projection $E\la E/\!\sim$ is
continuous (see for instance \cite{Bou}). The equivalence relation
$\sim$ is said to be {\rm closed} if its graph $\{(x,y)\in E\times
E\;:\;x\sim y\}$ is a closed subset of $E\times E$. We use several
times the following simple fact: if $E$ is a compact metric space and
$\sim$ is closed, then the quotient space $E/\!\sim$ is a Hausdorff
space, and is therefore compact, as the image of $E$ under the
canonical projection.

Let $\sphere^1$ be the unit circle in the complex plane $\CC$. If
$a,b\in\sphere^1$ and $a\neq b$, we denote by $[a,b]$ the closed arc in
$\sphere^1$ going from $a$ to $b$ in the counterclockwise order.
Similarly, $]a,b[$ denotes the corresponding open arc. By convention,
$[a,a]=\{a\}$ and $]a,a[\;=\emptyset$.

Let $g:\sphere^1\la \RR$ be a continuous function. For every
$a,b\in\sphere^1$, we set
$$
m_g(a,b)=
\max\left\{\min_{c\in[a,b]} g(c), \min_{c\in[b,a]} g(c)\right\}\;,
$$
and
$$
d_g(a,b)= g(a)+g(b)-2 m_g(a,b)\;.
$$
Note that $d_g(a,b)=0$ if and only if $g(a)=g(b)=m_g(a,b)$. We define
a closed equivalence relation $\sim_g$ on $\sphere^1$ by setting $a\sim_g
b$ if and only if $d_g(a,b)=0$.

Then $d_g$ induces a metric, still denoted by $d_g$, on the quotient
space ${\cal T}_g=\sphere^1/\!\!\sim_g$. Furthermore, ${\cal T}_g$
equipped with this metric is a compact $\RR$-tree. See Theorem 2.1 in
\cite{DuLG}, which deals with a slightly different but equivalent
setting. It is also easy to verify that the topology of the metric
space $({\cal T}_g,d_g)$ coincides with the quotient topology. Indeed,
the canonical projection $\sphere^1\la {\cal T}_g$ is continuous when
${\cal T}_g$ is equipped with the metric $d_g$, hence induces a
continuous bijection from ${\cal T}_g$ endowed with the quotient
topology, onto ${\cal T}_g$ endowed with the topology induced by the
metric $d_g$.  As $\t_g$ is compact for both topologies, the desired
result follows.

 From now on, we make the following additional assumption on $g$.

\begin{center}
$(H_g)$ ~~~~~~Local minima of $g$ are distinct.
\end{center}

\noindent %
This means that if $]a,b[$ and $]c,d[$ are two disjoint open arcs in
$\sphere^1$, and if the lower bound of the values of $g$ over $]a,b[$,
respectively over $]c,d[$, is attained at a point of $]a,b[$,
respectively at a point of $]c,d[$, then
$$
\min_{x\in[a,b]} g(x)\neq \min_{x\in[c,d]} g(x)\;.
$$

Let $\DD$ be the open unit disk in $\CC$ and let $\overline\DD
=\DD\cup\sphere^1$ be the closed disk. We equip $\DD$ with the usual
hyperbolic metric and for every $a,b\in\sphere^1$ with $a\neq b$, we
denote by $ab$ the (hyperbolic) geodesic line joining $a$ to $b$ in
$\DD$.  We also denote by $\overline{ab}$ the union of $ab$ and of the
points $a$ and $b$. By convention, $aa=\emptyset$ and $\overline{aa}
=\{a\}$.  We then let $L_g$ be the union of the geodesic lines $ab$ for
all pairs $\{a,b\}$ of distinct points of $\sphere^1$ such that $a\sim_g
b$.

Recall that a {\it (hyperbolic) geodesic lamination} in $\DD$ is a
closed subset of $\DD$ which is the union of a collection of pairwise
disjoint geodesic lines. A geodesic lamination is said to be {\it
  maximal} if it is maximal for the inclusion relation.  As a general
reference about geodesic laminations, we will use \cite{Bon} and the
references therein.

\begin{proposition} \label{prop:mexgeodlam} %
Under Assumption $(H_g)$, the set $L_g$ is a maximal geodesic
lamination of the hyperbolic disk
$\DD$.
\end{proposition}

\dem %
An elementary argument shows that, under assumption $(H_g)$,
equivalence classes for $\sim_g$ can have at most three points. Then,
let $\{a,b\}$ and $\{c,d\}$ be two pairs of distinct points in $\sphere^1$
such that $a\sim_g b$ and $c\sim_g d$. We claim that either the open
arcs $]a,b[$ and $]c,d[$ are disjoint, or one of them is contained in
the other one. Indeed, if this were not the case, then it would follow
from the definition of $d_g$ that the four points $a,b,c,d$ are
distinct and equivalent for $\sim_g$, which contradicts the first
observation of the proof. We conclude that the geodesic lines $ab$ and
$cd$ are disjoint, or coincide if $\{a,b\}=\{c,d\}$. Hence $L_g$ is a
disjoint union of geodesic lines.

As the equivalence relation $\sim_g$ is closed, its graph is compact
in $\sphere^1\times\sphere^1$. It immediately follows that $L_g$ is a closed
subset of the hyperbolic disk.

It remains to verify that $L_g$ is maximal. To this end, we argue by
contradiction. Let $a$ and $b$ be two distinct points in $\sphere^1$, and
suppose that $ab$ does not intersect $L_g$. Without loss of generality,
we may assume that
$$
\min_{x\in[a,b]} g(x)\geq \min_{x\in[b,a]} g(x)\;.
$$
If
$$
g(a)> \min_{x\in[a,b]} g(x)\;,
$$
an elementary argument shows that we can find two distinct points
$c\in\;]a,b[$ and $d\in\;]b,a[$ such that
$$
g(a)> g(c)=g(d)=\min_{x\in[c,d]} g(x)>\min_{x\in[a,b]} g(x)\;.
$$
But then $c\sim_g d$, and the geodesic line $cd$ intersects $ab$,
which contradicts our initial assumption that $ab$ does not intersect
$L_g$. We conclude that
$$
g(a)= \min_{x\in[a,b]} g(x)\;,
$$
and similarly we have
$$
g(b)= \min_{x\in[a,b]} g(x)\;.
$$
It follows that $a\sim_g b$, which is again a contradiction. \cq

\medskip %
Since $L_g$ is a maximal geodesic lamination, we know (see for
instance \cite{Bon}) that every connected component of $\DD\setminus
L_g$ is an ideal hyperbolic triangle. Clearly, these connected
components are in one-to-one correspondence with triples $\{a,b,c\}$
of distinct points in $\sphere^1$ such that $a\sim_g b\sim_g c$.

We can extend the equivalence relation $\sim_g$ to $\overline\DD$ as
follows. If $x,y\in\overline\DD$ and $x\neq y$, we put $x\sim_g y$ if
and only if $x$ and $y$ belong to the same arc $\overline{ab}$ with
$a\sim_g b$, or if $x$ and $y$ belong to the closure of the same ideal
geodesic triangle which is a connected component of $\DD\setminus
L_g$. In order to verify that this extension is still an equivalence
relation, we observe that a given geodesic line $ab$ cannot be
contained in the boundary of two distinct components of $\DD\setminus
L_g$. This again follows from the fact that equivalence classes for
$\sim_g$ contain at most three points of $\sphere^1$. For the extended
equivalence relation, equivalence classes are of three possible types,
either singletons $\{a\}$ for certain values of $a\in \sphere^1$, or arcs
$\overline{ab}$ for $a,b\in\sphere^1$, $a\neq b$ and $a\sim_g b$, or
closures of ideal hyperbolic triangles with ends $a,b,c$ such that
$a\sim_g b\sim_g c$.

By the preceding remarks, the inclusion map $\sphere^1\la \overline\DD$
induces a bijection $\sphere^1/\!\sim_g\;\la \overline\DD/\!\sim_g$, and
we use this bijection to identify these two sets. Note that this
identification is also an homeomorphism.  Indeed, the inclusion map
$\sphere^1\la \overline\DD$ is continuous and both $\sphere^1/\!\sim_g$ and
$\overline\DD/\!\sim_g$ are compact (note that the equivalence
relation $\sim_g$ on $\overline\DD$ is also closed).

The following two propositions are not used in the proofs of our main
results.  Still they contain useful information and answer basic
questions about the geodesic lamination $L_g$. We refer for instance
to \cite[page 12]{Bon} for the definition of a transverse measure on a
geodesic lamination, and to \cite[page 84]{LP} for the definition of
its space of leaves made Hausdorff.

\begin{proposition}\label{prop:lamigeod} %
Under Assumption $(H_g)$, the geodesic lamination $L_g$ carries a
natural transverse measure $\mu$, whose support is $L_g$, such that
the space of leaves made Hausdorff of $(L_g,\mu)$ is an $\RR$-tree
whose completion is isometric to the $\RR$-tree $(\t_g,d_g)$.
Furthermore, if the times of local minima of $g$ are dense in $\sphere^1$,
   then $L_g$ has empty interior.
\end{proposition}

\dem %
Let $\pi:\DD\la\overline\DD/\!\sim_g$ be the composition of the
inclusion map $\DD\la\overline\DD$ with the canonical projection
$\overline\DD\la\overline\DD/\!\sim_g$.  Consider in $\DD$ a non
trivial (hyperbolic) geodesic segment $[u,v]$, with $u,v\in\DD$, and
assume that this segment is transverse to $L_g$.  As a geodesic line
in $\DD$, that does not contain $[u,v]$, cuts (transversely) $[u,v]$
at one point at most, the restriction of the map $\pi$ to $L_g\cap
[u,v]$ is continuous and injective, except that the endpoints of a
connected component of $[u,v]\setminus L_g$ are mapped to the same
point. In particular, the image of this restriction is the geodesic
segment in $\t_g$ between $\pi(u)$ and $\pi(v)$. Denote by $\lambda$
the Lebesgue measure on this segment, which is isometrically
identified with an interval of the real line. Since $\lambda$ has no
atom, there exists a unique finite measure $\mu_{[u,v]}$ on $[u,v]$,
which is supported on $L_g\cap [u,v]$, such that the image measure of
$\mu_{[u,v]}$ under $\pi$ is $\lambda$. As the support of $\lambda$ is
$[\pi(u),\pi(v)]$, it follows that the support of $\mu_{[u,v]}$ is
exactly $L_g\cap [u,v]$. By construction, it is easy to check that the
transverse measure $\mu= (\mu_{[u,v]} )_{[u,v]}$ is invariant by
holonomy along the leaves of $L_g$. Hence $(L_g,\mu)$ is a
transversely measured geodesic lamination of $\DD$ (see
\cite{FLP,Bon}).

Now consider the pseudo-distance $\wt d$ on $\DD$, where $\wt d(u,v)$
is defined as the lower bound over all piecewise transverse paths
$\gamma$ from $u$ to $v$ of the total mass placed on $\gamma$ by the
transverse measure $\mu$.  Then the leaf space made Hausdorff
$T_{L_g,\mu}$ of $(L_g,\mu)$ is the quotient metric space of $(\DD,\wt
d)$ (obtained by identifying $u$ and $v$ if and only if $\wt
d(u,v)=0$), which is an $\RR$-tree (see \cite{MO,GS}).

Note that for every $u$ and $v$ in $\DD$, if the geodesic segment
$[u,v]$ is transverse to $L_g$, then $\wt d(u,v)=d_g(\pi(u),\pi(v))$,
as any piecewise transverse path from $u$ to $v$ has transverse
measure at least the transverse measure of $[u,v]$, by standard
arguments.  Hence the map $\pi$ induces an isometric embedding from
$T_{L_g,\mu}$ into $(\t_g,d_g)$. As $\DD$ is dense in $\overline\DD$,
the image of this embedding is dense. As $(\t_g,d_g)$ is compact, it
is hence (isometric to) the completion of $T_{L_g,\mu}$.

Suppose that the times of local minima of $g$ are dense. Let
$a,b\in\sphere^1$ be such that $a\sim_g b$. Without loss of generality,
assume that $\min_{x\in[a,b]} g(x)\geq \min_{x\in[b,a]} g(x)$. For
every $c$ in $]a,b[$, if $d$ is a local minimum of $g$ that belongs to
$]a,c[$, then there exists $e\in\;[a,d]\setminus\{d\}$ such that
$d\sim_g e$. As distinct geodesic lines in $L_g$ are disjoint, no
point of the geodesic line $ab$ can be in the interior of $L_g$.  \cq

\medskip %
If $A$ is a subset of the closed disk $\ov\DD$ equipped with the usual
Euclidean distance, we denote by ${\rm dim}(A)$ the Hausdorff
dimension of $A$, and by $\underline{\rm dim}_M(A)$ the lower
Minkowski dimension of $A$ (also called the lower box-counting
dimension, see for instance \cite[page 77]{Mat}). Recall that ${\rm
  dim}(A) \leq \underline{\rm dim}_M(A)$.

Let $A_g$ denote the set of all $x\in \sphere^1$ such that the equivalence
class of $x$ under $\sim_g$ is not a singleton. We also let $\cal J$
be the countable set of all (ordered) pairs $(I,J)$ where $I$ and $J$
are two disjoint closed subarcs of $\sphere^1$ with nonempty interior and
rational endpoints. If $(I,J)\in{\cal J}$, we denote by $A^{(I,J)}_g$
the set of all $x\in I$ such that $x\sim_g y$ for some $y\in J$.
Plainly,
\begin{equation}
\label{Haus-decomp}
A_g=\bigcup_{(I,J)\in {\cal J}} A_g^{(I,J)}\;.
\end{equation}

\begin{proposition} \label{Hausdorff-lami} %
{\rm(i)} We have 
$$
{\rm dim}(L_g)\geq 1+{\rm dim}(A_g)\;.
$$ 
{\rm(ii)} Assume that $\underline{\rm dim}_M(A_g^{(I,J)}\cup
A_g^{(J,I)})\leq {\rm dim}(A_g)$ for every $(I,J)\in {\cal J}$. Then,
$$
{\rm dim}(L_g)= 1+{\rm dim}(A_g)\;.
$$
\end{proposition}

\dem %
(i) We assume that ${\rm dim}(A_g)>0$, because otherwise the result is
easy. Let $\alpha\in\,]0,{\rm dim}(A_g)[$. By (\ref{Haus-decomp}), we can find a pair
$(I,J)\in{\cal J}$ such that ${\rm dim}(A_g^{(I,J)})>\alpha$. Since
$A_g^{(I,J)}$ is a compact subset of $\sphere^1$, Frostman's lemma
\cite[page 112]{Mat} yields the existence of a nontrivial finite
Borel measure $\nu$ supported on $A_g^{(I,J)}$ such that
$$
\nu(B(x,r))\leq r^\alpha\eqno{(2)}
$$ 
for every $r>0$ and $x\in \sphere^1$. Here $B(x,r)$ denotes the
(Euclidean) disk of radius $r$ centered at $x$. Let $\widetilde
A_g^{(I,J)}$ be the subset of $A_g^{(I,J)}$ consisting of points $x$
such that the equivalence class of $x$ under $\sim_g$ contains exactly
two points, and for every $x\in \widetilde A_g^{(I,J)}$, let $s_g(x)$
be the unique element of $J$ such that $x\sim_g s_g(x)$. Notice that
$A_g^{(I,J)}\backslash \widetilde A_g^{(I,J)}$ is countable, and so
$\nu$ is supported on $\widetilde A_g^{(I,J)}$.  For every $x\in
\widetilde A_g^{(I,J)}$, let $\lambda_x$ denote the one-dimensional
Hausdorff measure on the arc $xs_g(x)$ (equipped with the Euclidean
distance). Define a finite Borel measure $\Lambda$ by setting for
every Borel subset $B$ of the plane
$$
\Lambda(B)=\int\nu(dx)\int \lambda_x(dz)\,{\bf 1}_B(z).
$$
By construction, $\Lambda$ is supported on $L_g$. Then fix
$R\in\,]0,1[$ such that $\Lambda(B(0,R))>0$. Let $z_0\in L_g$ be such
that $|z_0|\leq R$, and choose $x_0,y_0\in A_g$ such that $z_0\in
x_0y_0$. Let $\varepsilon\in\,]0,1]$. A simple geometric argument
shows that the conditions $x\in \widetilde A_g^{(I,J)}$ and
$xs_g(x)\cap B(z_0,\varepsilon)\not =\emptyset$ imply $|x-x_0|\leq
C\varepsilon$, where the constant $C$ only depends on $R$.  Hence,
using (2),
$$
\Lambda(B(z_0,\varepsilon))=\int_{\{|x-x_0|\leq C\varepsilon\}}
\nu(dx)\,\lambda_x(B(z_0,\varepsilon))
\leq C'\varepsilon^{1+\alpha}
$$
where the constant $C'$ does not depend on $\varepsilon$ nor on $z_0$.
Frostman's lemma now gives ${\rm dim}(L_g)\geq 1+\alpha$ as desired.

\noindent(ii) We now prove that ${\rm dim}(L_g)\leq 1+{\rm dim}(A_g)$
under the assumption in (ii). For $(I,J)\in {\cal J}$, let
$F^{(I,J)}_g$ be the union of all geodesic lines $xy$ for $x\in I$,
$y\in J$ and $x\sim_gy$. It is enough to prove that ${\rm
  dim}(F^{(I,J)}_g)\leq 1+{\rm dim}(A_g)$ for a fixed choice of
$(I,J)\in {\cal J}$.

Let $\beta>{\rm dim}( A_g)$. By the assumption in (ii), we can find a
sequence $\varepsilon_1>\varepsilon_2>\cdots>0$ decreasing to $0$,
such that the following holds for every $\varepsilon$ belonging to
this sequence.  There exist a positive integer $M(\varepsilon)\leq
\varepsilon^{-\beta}$ and $M(\varepsilon)$ disjoint subarcs
$I_1,I_2,\ldots,I_{M(\varepsilon)}$ of $\sphere^1$, with length less than
$\varepsilon$, such that $A^{(I,J)}_g$ is contained in
$I_1\cup\cdots\cup I_{M(\varepsilon)}$. Similarly, we can find a
positive integer $N(\varepsilon)\leq \varepsilon^{-\beta}$ and
$N(\varepsilon)$ disjoint subarcs $J_1,J_2,\ldots,J_{N(\varepsilon)}$,
with length less than $\varepsilon$, such that $A^{(J,I)}_g$ is
contained in $J_1\cup\cdots\cup J_{N(\varepsilon)}$. Then, let $H$ be
the set of all pairs $(i,j)\in\{1,\ldots,M(\varepsilon)\}\times
\{1,\ldots,N(\varepsilon)\}$ such that there exists a geodesic line
$xy\subset L_g$ with $x\in I_i$ and $y\in J_j$. Because geodesic lines
in $L_g$ are not allowed to cross, a simple argument shows that
$\#(H)\leq M(\varepsilon)+N(\varepsilon)\leq 2\varepsilon^{-\beta}$.
It easily follows that the two-dimensional Lebesgue measure of the
Euclidean tubular neighborhood of $F^{(I,J)}_g$ with radius
$\varepsilon$ is bounded above by $C\varepsilon^{1-\beta}$, where the
constant $C$ does not depend on $\varepsilon$ in our sequence.  This
implies \cite[page 79]{Mat} that $\underline{\rm dim}_M(F^{(I,J)}_g)
\leq 1+\beta$, and a fortiori ${\rm dim}(F^{(I,J)}_g)\leq 1+\beta$.
\cq

\bigskip %
We now come to the main result of this section.  We let $h:\sphere^1
\la\RR$ be another continuous function. We again assume that local
minima of $h$ are distinct, i.e.~that $(H_h)$ holds.  Furthermore, we
assume that the following condition holds.

\begin{center}
    $(H'_{g,h})$ ~Let $a,b,c$ be three points in $\sphere^1$ such that
    $a\sim_g b$ and $a\sim_h c$. \\ Then $a=b$ or $a=c$.
\end{center}

\noindent %
In other words, if the equivalence class of $a\in\sphere^1$ with respect
to $\sim_g$ is not a singleton, then its equivalence class with
respect to $\sim_h$ must be a singleton.

We can define an equivalence relation, which we still denote by
$\sim_h$, on the quotient $\t_g=\sphere^1/\!\sim_g$ by declaring for
$\alpha,\beta\in\t_g$ that $\alpha\sim_h \beta$ if and only if there
exists a representative $a$ of $\alpha$ in $\sphere^1$, respectively a
representative $b$ of $\beta$ in  $\sphere^1$, such that $a\sim_h b$.
Note that our assumption $(H'_{g,h})$ is used to verify that this
prescription defines an equivalence relation on $\t_g$.

\begin{proposition} \label{prop:quotspher} %
Under Assumptions $(H_g)$, $(H_h)$ and $(H'_{g,h})$,
the quotient space $\t_g/\!\sim_h$ is homeo\-morphic to the sphere
$\sphere^2$.
\end{proposition}

\dem %
We embed the complex plane into $\RR^3$ by identifying it with the
horizontal plane $\{x_3=0\}$. We write $H_+=\{x=(x_1,x_2,x_3)\in \sphere^2
\;:\; x_3\geq 0\}$ for the (closed) upper hemisphere, and similarly
$H_-=\{x=(x_1,x_2,x_3)\in\sphere^2\;:\; x_3\leq 0\}$ for the lower
hemisphere.  We can use the stereographic projection from the South
pole to identify (topologically) $H_+$ with the closed unit disk
$\overline\DD$.  Thanks to this identification, we may define the
equivalence relation $\sim_g$ on $H_+$, and by previous observations,
the quotient space $H_+/\!\sim_g$ is homeomorphic to $\t_g$.
Similarly, we can use the stereographic projection from the North pole
to identify $H_-$ with $\overline\DD$, and then define the equivalence
relation $\sim_h$ on $H_-$.

Let $\sim$ be the equivalence relation on $\sphere^2$ whose graph is the
union of the graphs of $\sim_g$ and $\sim_h$ viewed as equivalence
relations on $H_+$ and $H_-$ respectively. Note that $(H'_{g,h})$ is
used to verify that $\sim$ is an equivalence relation on $\sphere^2$. Any
equivalence class for $\sim$ is an equivalence class for $\sim_g$, or
an equivalence class for $\sim_h$. It may be both if and only if
it is a singleton. As a consequence, any equivalence class of $\sim$
is a compact path-connected subset of $\sphere^2$ whose complement is also
connected. Furthermore, as $\sim_g$ and $\sim_h$ are both closed, it
follows that the equivalence relation $\sim$ is closed.

At this point, we use the following theorem of Moore \cite[page
416]{Moo} (see also \cite{Thu} for a previous application of this
theorem, and in particular Figure 10, page 376 of \cite{Thu}).

\begin{theorem}[Moore] %
\label{moore-theorem}
Let $\sim$ be a closed equivalence relation on $\sphere^2$. Assume that
every equivalence class of $\sim$ is a compact path-connected subset
of the sphere whose complement is connected. Then the quotient space
$\sphere^2/\!\sim$ is homeomorphic to $\sphere^2$.
\end{theorem}

Clearly Moore's theorem applies to our setting, and we get that the
quotient $\sphere^2/\!\sim$ is homeomorphic to $\sphere^2$.

To complete the proof of Proposition \ref{prop:quotspher}, it remains
to verify that $\t_g/\!\sim_h$ is homeomorphic to $\sphere^2/\!\sim$. We
first observe that $\t_g/\!\sim_h$ is compact. Indeed, $\sim_h$ viewed
as an equivalence relation on $\t_g$ is closed, as $\sim_h$ is closed
on $\sphere^1$ and the canonical projection $\sphere^1\la \t_g$ is a closed
map.  Then, by composing the inclusion map $\sphere^1\la \sphere^2$ with the
projection $\sphere^2\la \sphere^2/\!\sim$, we get a continuous mapping, which
factorizes through the equivalence relation $\sim_g$ and thus yields a
continuous mapping from $\t_g$ onto $\sphere^2/\!\sim$. Again, this
mapping factorizes through the equivalence relation $\sim_h$ and we
obtain that the canonical bijection from $\t_g/\!\sim_h$ onto
$\sphere^2/\!\sim$ is continuous.  Since both $\t_g/\!\sim_h$ and
$\sphere^2/\!\sim$ are compact, this bijection is a homeomorphism.
\cq

\medskip 
\rem %
Assumption $(H'_{g,h})$ in Proposition 2.4 can be weakened. The
application of Moore's theorem is possible under less stringent
assumptions.

\section{Proof of the main result}
\label{sec:proofmain}

In this section we prove Theorem \ref{theo:intromain} and Corollary
\ref{coro:introtype}.  On a given probability space, we consider a
normalized Brownian excursion $(\eg_t)_{0\leq t\leq 1}$ and a process
$(Z_t)_{0\leq t\leq 1}$ which is distributed as the head of the
one-dimensional Brownian snake driven by $\eg$. This means that the
process $Z$ has continuous sample paths and that, conditionally given
$\eg$, it is a centered real-valued Gaussian process with
(conditional) covariance function
\begin{equation}
\label{cov}
E[Z_sZ_t\mid\eg]=\min_{s\wedge t\leq u\leq s\vee t} \eg_u
\end{equation}
for every $s,t\in [0,1]$. See Section 4 below for more information
about $Z$ and the Brownian snake. Notice that $\eg_0=\eg_1=0$
and $Z_0=Z_1=0$, a.s.

We also need to introduce the pair $(\eg,Z)$ ``re-rooted at the
minimal spatial position''. Set
$$
\un Z=\min_{0\leq s\leq 1} Z_s
$$
and let $s_*$ be the almost surely unique time such that $Z_{s_*}=\un
Z$ (the uniqueness of $s_*$ follows from Proposition 2.5 in
\cite{LGW}, and is also a consequence of Lemma \ref{localminima}
below).  For every $s,t\in[0,1]$, set $s\oplus t=s+t$ if $s+t\leq 1$
and $s\oplus t=s+t-1$ if $s+t>1$. Then, for every $s\in[0,1]$, define
\begin{description}
\item{$\bullet$}
$\displaystyle{\ov\eg_t=\eg_{s_*}+\eg_{s_*\oplus
t}-2\min_{s_*\wedge(s_*\oplus s)
\leq r \leq s_* \vee (s_*\oplus s)} \eg_r}$;
\item{$\bullet$} $\ov Z_t=Z_{s_*\oplus t} -Z_{s_*}$.
\end{description}
Note again that $\ov\eg_0=\ov\eg_1=0$ and $\ov Z_0=\ov Z_1=0$ a.s. The
pair $(\ov \eg,\ov Z)$ can be interpreted as the pair $(\eg,Z)$
conditioned on the event $\{\un Z=0\}$ (see \cite{LGW}).

In view of applying the results of Section 2, it will be convenient to
view the random functions $\eg,Z,\ov \eg$ and $\ov Z$ as parametrized
by the circle $\sphere^1$ rather than by the interval $[0,1]$. This is of
course easily achieved by setting, for instance,
$$
\eg(e^{2i\pi r})=\eg_r\ ,\quad r\in [0,1].
$$

Then Assumption ${(H_\eg)}$ holds a.s. This follows from the
well-known analogous result for linear Brownian motion, which is a
very easy application of the Markov property.

\begin{lemma}
\label{localminima}
Assumption ${(H_{Z})}$ holds almost surely. In other words, local
minima of $Z$ are distinct, with probability one.
\end{lemma}

\begin{lemma}
\label{keylemma}
Asumption ${(H'_{\eg,{Z}})}$ holds almost surely. In other words,
almost surely for every $a,b,c\in \sphere^1$, the conditions $a\sim_\eg b$
and $a\sim_Z c$ imply that $a=b$ or $a=c$.
\end{lemma}

We postpone the proof of these two lemmas to Section 4. Thanks to
these lemmas, we can apply the results of Section 2 to the pair
$(\eg,Z)$. In particular, we can consider the quotient space
$\t_\eg\,/\!\sim_Z$ and we know from Proposition \ref{prop:quotspher}
that this quotient space is almost surely homeomorphic to the sphere
$\sphere^2$.

In order to complete the proof of Theorem \ref{theo:intromain}, it
will be sufficient to verify that the (random) metric space that
appears in \cite{L} as the weak limit of rescaled random maps is a.s.
homeomorphic to $\t_\eg\,/\!\sim_Z$. We first need to recall the
topological description of this limiting random metric space that is
given in \cite{L}.

We start by observing that outside a set of probability zero, the
value of $Z_a$ for $a\in \sphere^1$, respectively the value of $\ov Z_a$,
only depends on the equivalence class of $a$ in $\t_\eg$, respectively
in $\t_{\ov\eg}$\,: This essentially follows from the form of the
covariance function in (\ref{cov}), see Section 2.4 in \cite{L}.
Thanks to this observation, we may and will sometimes view $Z$,
respectively $\ov Z$, as parametrized by $\t_\eg$, respectively
$\t_{\ov\eg}$.

Let us denote by $p_\eg$, respectively $p_{\ov\eg}$, the canonical
projection from $\sphere^1$ onto $\t_\eg$, respectively onto
$\t_{\ov\eg}$. If $\alpha,\beta\in \t_\eg$, we denote by
$[\alpha,\beta]$ the image under $p_\eg$ of the smallest arc $[a,b]$
in $\sphere^1$ such that $p_\eg(\alpha)=a$ and $p_\eg(\beta)=b$. We
similarly define $[\alpha,\beta]$ when $a,\beta\in \t_{\ov\eg}$. Then,
for every $\alpha,\beta\in \t_{\ov\eg}$, we set
$$
D^\circ (\alpha,\beta)=\ov Z_\alpha+\ov Z_\beta
-2\max\Big(\min_{\gamma\in[\alpha,\beta]} \ov Z_\gamma,
\min_{\gamma\in[\beta,\alpha]} \ov Z_\gamma\Big)
$$
and
$$
D^*(\alpha,\beta)=\inf\Big\{\sum_{i=1}^p D^\circ
(\alpha_{i-1},\alpha_i)\Big\}
$$
where the last lower bound is over all choices of the integer $p\geq
1$ and of the finite sequence $\alpha_0,\alpha_1,\ldots,\alpha_p$ in
$\t_{\ov\eg}$, such that $\alpha_0=\alpha$ and $\alpha_p=\beta$. We
set $\alpha\approx\beta$ if and only if $D^*(\alpha,\beta)=0$. By
Theorem 3.4 in \cite{L}, this is also equivalent to the condition
$D^\circ(\alpha,\beta)=0$, for every $\alpha,\beta\in\t_{\ov\eg}$, almost
surely.

Recall the notation introduced in Section 1. According to the same
theorem of \cite{L} and Remark (a) following it, any weak limit of the
sequence $({\bf m}_n,n^{-1/4}d_n)$ is a.s. homeomorphic to the
quotient space $\t_{\ov\eg}\,/\!\approx$ equipped with the metric
induced by $D^*$, which is still denoted by $D^*$. Thus Theorem
\ref{theo:intromain} follows from the next proposition.

\begin{proposition}
\label{techmail}
The metric space $(\t_{\ov\eg}\,/\!\approx,D^*)$ is almost surely
homeomorphic to the quotient space $\t_\eg\,/\!\sim_Z$.
\end{proposition}

\dem
We first construct a (canonical) bijection between
$\t_{\ov\eg}\,/\!\approx$ and $\t_\eg\,/\!\sim_Z$ and then verify that
this bijection is a homeomorphism.  Let $\rho:\sphere^1\la \sphere^1$ be the
rotation with angle $2\pi s_*$.  According to the re-rooting lemma
(Lemma 2.2 in \cite{DuLG}), $\rho$ induces an isometry $R$ from
$(\t_{\ov\eg},d_{\ov\eg})$ onto $(\t_{\eg},d_\eg)$. Furthermore, for
every $\alpha\in \t_\eg$,
\begin{equation}
\label
{reroot-spatial}
\ov Z_\alpha = Z_{R(\alpha)}-\un Z.
\end{equation}
Now recall that a.s. for every $\alpha,\beta\in \t_{\ov\eg}$, the
relation $\alpha\approx\beta$ holds if and only if
$D^\circ(\alpha,\beta) =0$, or equivalently
$$
\ov Z_\alpha =\ov Z_\beta =\max\Big(\min_{\gamma\in[\alpha,\beta]}
\ov Z_\gamma,
\min_{\gamma\in[\beta,\alpha]} \ov Z_\gamma\Big)\;.
$$
 From our definitions and the identity (\ref{reroot-spatial}), this is
immediately seen to be equivalent to
$$
Z_{R(\alpha)} =Z_{R(\beta)}
=\max\Big(\min_{\gamma\in[R(\alpha),R(\beta)]} Z_\gamma,
\min_{\gamma\in[R(\beta),R(\alpha)]} Z_\gamma\Big)\;,
$$
that is to $R(\alpha)\sim_Z R(\beta)$.

Thus $R$ induces a bijection, which we denote by $\wt R$, from
$\t_{\ov\eg}\,/\!\approx$ onto $\t_\eg\,/\!\sim_Z$. To prove that $\wt
R$ is a homeomorphism, it is enough to verify that $\wt R^{-1}$ is
continuous, since both $\t_{\ov\eg}\,/\!\approx$ (equipped with the
metric $D^*$) and $\t_\eg\,/\!\sim_Z$ are compact. The canonical
projection from $\t_{\ov\eg}$ onto $(\t_{\ov\eg}\,/\!\approx,D^*)$ is
continuous: Using the continuity of the mapping $\sphere^1\ni a\to \ov
Z_a$, a direct inspection of the definition of $D^\circ$ shows that if
$\alpha_n$ tends to $\alpha$ in $\t_{\ov{\bf e}}$ then
$D^\circ(\alpha_n,\alpha)$ tends to $0$ and a fortiori
$D^*(\alpha_n,\alpha)$ tends to $0$ as $n\to\infty$. By composing the
isometry $R^{-1}$ from $(\t_\eg,d_{\eg})$ onto
$(\t_{\ov\eg},d_{\ov\eg})$ with the previous projection, we get a
continuous mapping from $(\t_\eg,d_{\eg})$ onto
$(\t_{\ov\eg}\,/\!\approx,D^*)$, which in turn induces a continuous
mapping from the space $\t_\eg\,/\!\sim_Z$, equipped with the quotient
topology, onto $(\t_{\ov\eg}\,/\!\approx,D^*)$. The latter mapping is
just $\wt R^{-1}$, and so we have obtained that $\wt R^{-1}$ is
continuous, which completes the proof.  \cq

\medskip
\noindent{\bf Proof of Corollary \ref{coro:introtype}.} For every
integer $n\geq 2$, denote by $A_n$ the event consisting of all
$\omega$'s in our underlying probability space such that there exists
an injective cycle of the map $M_n(\omega)$ satisfying the properties
stated in the corollary. We argue by contradiction, assuming that
$P(A_n)$ does not converge to $0$. Then we can find $\eta>0$ and a
sequence $(n_k)$ converging to $+\infty$ such that $P(A_{n_k})\geq
\eta$ for every $k$. From now on, we restrict our attention to values
of $n$ belonging to this sequence. By extracting another subsequence
if necessary, we can also assume that $({\bf m}_n,n^{-1/4}d_n)$
converges in distribution along this sequence. The convergence in
distribution can be replaced by an almost sure convergence thanks to
the Skorokhod representation theorem.  Thus we have almost surely
\begin{equation}
\label{convGH}
({\bf m}_n,n^{-1/4}d_n)\longrightarrow ({\bf m}_\infty,d_\infty)
\end{equation}
as $n\to\infty$, in the sense of the Hausdorff-Gromov distance. By
Theorem \ref{theo:intromain}, $({\bf m}_\infty,d_\infty)$ is almost
surely homeomorphic to the sphere $\sphere^2$.

From now on, we argue with a fixed value of $\omega$ in our
probability space, such that $\omega\in \limsup A_n$ (this event has
probability greater than $\eta$ by the above), the convergence
(\ref{convGH}) holds and $({\bf m}_\infty,d_\infty)$ is homeomorphic
to $\sphere^2$. Let us show that this leads to a contradiction. By the
definition of the events $A_n$, we can find a subsequence (depending
on $\omega$) such that for every $n$ belonging to this subsequence,
there exists an injective cycle $C_n$ of the map $M_n$, with length
$\ell(C_n)\leq \theta(n)$ and two vertices $a_n,b_n\in {\bf m}_n$,
which are separated by the cycle $C_n$ (in the sense that every
(continuous) path from $a_n$ to $b_n$ has to cross $C_n$) and such
that $\min\{d_n(a_n,C_n), d_n(b_n,C_n)\} > \delta n^{1/4}$. Here
$d_n(a_n,C_n)$ denotes as usual the minimal distance between $a_n$ and
a vertex of $C_n$.

Say that a map $\varphi$ from a metric space $(E,d)$ into another
metric space $(E',d')$ is an $\varepsilon$-isometry if
$|d'(\varphi(x),\varphi(y))-d(x,y)|\leq \varepsilon$ for every $x,y\in
E$. From the convergence (3) and the definition of the
Hausdorff-Gromov topology (see e.g.~\cite{BBI}), we can find a
sequence $\varepsilon_n\to 0$, and $\varepsilon_n$-isometries
$f_n:({\bf m}_n,n^{-1/4}d_n)\to({\bf m}_\infty,d_\infty)$ and
$g_n:({\bf m}_\infty, d_\infty)\to({\bf m}_n,n^{-1/4}d_n)$ such that
$n^{-1/4}d_n(g_n\circ f_n(x),x) \leq \epsilon_n$ for every $x$ in
${\bf m}_n$.  Let $a'_n=f_n(a_n), b'_n=f_n(b_n)$ and let $C'_n$ be the
image under $f_n$ of the vertex set of $C_n$. Note that the diameter
of $C'_n$ tends to $0$ by our assumption $\ell(C_n)\leq \theta(n)$.
Using the compactness of ${\bf m}_\infty$ and again extracting a
subsequence if necessary, we can assume that the points $a'_n,b'_n$
converge respectively to $a_\infty,b_\infty$ in ${\bf m}_\infty$ and
the (finite) sets $C'_n$ converge (for the Hausdorff distance) to a
singleton $\{c_\infty\}$, such that $\min\{ d_\infty
(a_\infty,c_\infty), d_\infty(b_\infty,c_\infty)\}\geq \delta$.  Since
the complement of a single point in the sphere $\sphere^2$ is path
connected, there exists a (continuous) path $\gamma:[0,1]\to {\bf
  m}_\infty$ from $a_\infty$ to $b_\infty$ avoiding $c_\infty$. Let
$\epsilon'= \min\{\varepsilon, d_\infty(c_\infty,\gamma)\}>0$, and let
$N\in\NN$ be large enough so that $d_\infty\big(\gamma(\frac{k}{N}),
\gamma(\frac{k+1}{N})\big) \leq \frac{\epsilon'}{3}$ for every integer
$k$ with $0\leq k\leq N-1$.  For $0\leq k\leq N-1$, define $x_{n,k}=
g_n(\gamma(\frac{k}{N}))$, and $x_{n,-1}=a_n$, $x_{n,N+1}=b_n$. Then
if $n$ is large enough, $(x_{n,k})_{-1\leq k\leq N+1}$ is a sequence
of points in ${\bf m}_n$ such that $d_n(x_{n,k},C_n)\geq
\frac{\epsilon'}{2}n^{1/4}$ and $d_n(x_{n,k},x_{n,k+1}) <
\frac{\epsilon'}{2}n^{1/4}$. Connecting $x_{n,k}$ and $x_{n,k+1}$ by a
geodesic path in the graph $M_n$, we get a path from $a_n$ to $b_n$ in
the map $M_n$ that avoids $C_n$, which is a contradiction.  \cq

\medskip
We conclude this section with an application of Proposition
\ref{Hausdorff-lami} to the Hausdorff dimension of the random
lamination $L_{\bf e}$. The result is already stated
in \cite{Al4}, but the proof there is only sketched.

\begin{proposition}
\label{Hausdorff-random}
We have ${\rm dim}(L_{\bf e})=3/2$ almost surely.
\end{proposition}

\dem %
Denote by $\overline{\rm dim}_M(A)$ the upper Minkowski dimension of a
set $A$, and recall that $\overline{\rm dim}_M(A\cup B) =
\max(\overline{\rm dim}_M(A),\overline{\rm dim}_M(B))$. Also recall
the notation introduced before Proposition \ref{Hausdorff-lami}.  It
is then enough to prove that ${\rm dim}(A_{\bf e})=1/2$ and
$\overline{\rm dim}_M(A^{(I,J)}_{\bf e})\leq 1/2$, for every $(I,J)\in
{\cal J}$, a.s. In this proof, it is more convenient to view ${\bf e}$
as parametrized by the time interval $[0,1]$. Recall that ${\bf
  e}_0={\bf e}_1=0$ and ${\bf e}_t>0$ for every $t\in\;]0,1[$, a.s.,
and set
$$
H_a:=\{t\in\;]a,1]:{\bf e}_t=\min_{a\leq r\leq t} {\bf e}_r\}\;.
$$
for every $a\in\;]0,1[$. Then $H_a\subset A_e$ for every
$a\in\;]0,1[$, and $A^{(I,J)}_{\bf e}\subset H_a$ whenever $(I,J)\in {\cal J}$
and $a\in]0,1[$ are such that $I=[u,v]$, $J=[u',v']$ and 
$0\leq u'<v'<a<u<v\leq 1$. Using the invariance of
the Brownian excursion under time reversal, we hence see that the
required properties follow from the identities
$$
{\rm dim}(H_a)=\overline{\rm dim}_M(H_a)=\frac{1}{2}\ ,\qquad
\hbox{for almost all }a\in\,]0,1[,\quad\hbox{a.s.}
$$
(almost all refers to Lebesgue measure on $]0,1[$). By a scaling
argument, it is enough to verify that a similar property holds under
the It\^o measure of Brownian excursions. Using the Markov property at
time $a>0$ under the It\^o measure, it then suffices to prove that the
following holds. If $(\beta_t)_{t\geq 0}$ is a standard linear
Brownian motion, and
$$
K_u:=\{t\in[0,u]:\beta_t=\min_{0\leq r\leq t} \beta_r\}\;,
$$
we have
$$
{\rm dim}(K_u)=\overline{\rm dim}_M(K_u)=\frac{1}{2}
$$
for every $u>0$, a.s. By a classical theorem of L\'evy, the random set
$K_u$ has the same distribution as $\{t\in[0,u]:\beta_t=0\}$. The
preceding claim now follows from standard results about the zero set
of linear Brownian motion.  
\cq

\section{Proof of the technical results}

In this section, we prove Lemma \ref{localminima} and Lemma
\ref{keylemma}.  In both these proofs, it is more convenient to view
the processes $\eg$ and $Z$ as parametrized by the interval $[0,1]$
rather than by the unit circle (see the beginning of Section 3).  We
will make extensive use of properties of the Brownian snake.  We start
with a brief discussion of this path-valued Markov process, referring
to \cite{Zu} for a more thorough presentation (see also Section 4 in
\cite{L}).

A finite path in $\RR$ is a continuous mapping $\w:[0,\zeta]\la \RR$,
where $\zeta=\zeta_{(\w)}$ is a nonnegative real number called the
lifetime of $\w$. The set $\W$ of all finite paths is a Polish space
when equipped with the distance
$$
d(\w,\w')=|\zeta_{(\w)}-\zeta_{(\w')}|+\max_{t\geq 0}|\w(t\wedge
\zeta_{(\w)})-\w'(t\wedge\zeta_{(\w')})|\;.
$$

Let $x\in\RR$. The one-dimensional Brownian snake with initial point
$x$ is a continuous strong Markov process taking values in the set
$\W_x=\{\w\in\W:\w(0)=x\}$. Thus for every $s\geq 0$, $W_s=
(W_s(t),0\leq t\leq \zeta_s)$ is a random continuous path in $\RR$,
with a (random) lifetime $\zeta_s=\zeta_{(W_s)}$ and such that
$W_s(0)=x$.  The behavior of the Brownian snake can be described
informally as follows. The lifetime $\zeta_s$ evolves like reflecting
Brownian motion in $\RR_+$, and when $\zeta_s$ decreases the path
$W_s$ is erased from its tip, whereas when $\zeta_s$ increases the
path $W_s$ is extended by adding little pieces of Brownian paths at
its tip.  We denote by $\wh W_s=W_s(\zeta_s)$ the terminal point (head
of the snake) of the path $W_s$.

Let us fix $\w\in\W_x$. We denote by $\PP_\w$ the law of the Brownian
snake started from $\w$. We also let $\NN_x$ be the excursion measure
of the Brownian snake away from the trivial path with lifetime $0$ in
$\W_x$. Both measures $\PP_\w$ and $\NN_x$ may be defined on the space
$C(\RR_+,\W)$ of all continuous functions from $\RR_+$ into $\W$.
Under $\PP_\w$, the lifetime process $(\zeta_s,s\geq0)$ evolves like
reflecting Brownian motion in $\RR_+$, whereas under $\NN_x$ it is
distributed according to the It\^o measure of positive excursions of
linear Brownian motion.  In particular the quantity
$\sigma:=\sup\{s\geq 0:\zeta_s>0\}$ is finite $\NN_x$ a.e., and is
called the duration of the excursion. We denote by $\NN^{(1)}_x$ the
It\^o measure conditioned on the event $\{\sigma=1\}$.  Then, the pair
$(\eg_s,Z_s)_{0\leq s\leq 1}$ of Section 3 has the same distribution
as $(\zeta_s,\wh W_s)_{0\leq s\leq 1}$ under $\NN^{(1)}_0$, which
explains why Lemma \ref{localminima} and Lemma \ref{keylemma} will be
reduced to statements about the Brownian snake.

For every $s,s'\geq 0$, set
$$
m(s,s')=\min_{s\wedge s'\leq r\leq s\vee s'} \zeta_r\;.
$$
We use several times the so-called {\it snake property}: $\PP_\w$
a.s., or $\NN_x$ a.e.  (or $\NN^{(1)}_x$ a.s.) for every $s,s'\geq 0$,
we have
$$
W_s(t)=W_{s'}(t)\ ,\ \hbox{for every }0\leq t\leq m(s,s')\;.
$$

We now state a lemma which plays an important role in our proofs.
Recall that we have fixed an element $\w$ of $\W_x$.  Under the
probability measure $\PP_\w$, we set
$$
T=\inf\{s\geq 0:\zeta_s=0\}
$$
and we denote by $]\alpha_i,\beta_i[$, $i\in I$, the connected
components of the open set
$$
\{s\in[0,T]:\zeta_s>m(0,s)\}\;.
$$
Then, for every $i\in I$, we define a random element $W^i$ of
$C(\RR_+,\W)$ by setting for every $s\geq 0$,
$$
W^i_s(t)=W_{(\alpha_i+s)\wedge\beta_i}(\zeta_{\alpha_i}+t)\;
,\qquad 0\leq
t\leq\zeta^i_s:=\zeta_{(\alpha_i+s)\wedge\beta_i}-\zeta_{\alpha_i}\;.
$$
 From the snake property, $W^i_s$ belongs to ${\cal W}_{
   \w(\zeta_{\alpha_i})}$. The following result is Lemma V.5 in
\cite{Zu}.

\begin{lemma}
\label{subtrees}
The point measure
$$
\sum_{i\in I}\delta_{(\zeta_{\alpha_i},W^i)}(dt\,d\omega)
$$
is under $\PP_\w$ a Poisson point measure on $\RR_+\times C(\RR_+,\W)$
with intensity
$$
2\;{\bf 1}_{[0,\zeta_{(\w)}]}(t)dt\;\NN_{\w(t)}(d\omega)\;.
$$
\end{lemma}

We will also use the explicit form of the law of the minimal value of
the Brownian snake under $\NN_x$: For every $x,y\in \RR$ with $y<x$,
\begin{equation}
\label{minisnake}
\NN_x\Big(\min_{s\geq 0} \wh W_s<y\Big)=\frac{3}{2}\,(x-y)^{-2}\;.
\end{equation}
See e.g.~Lemma 2.1 in \cite{LGW}.

\smallskip
\noi{\bf Proof of Lemma \ref{localminima}.} %
Using the fact that the distribution of $(\eg_s,Z_s)_{0\leq s\leq 1}$
is the same as that of $(\zeta_s,\wh W_s)_{0\leq s\leq 1}$ under
$\NN^{(1)}_0$, together with a simple scaling argument, it is enough
to prove that local minima of $\wh W_s$ over the time interval
$[0,\sigma]$ are distinct $\NN_0$ a.e. This is less easy than the
analogous result for linear Brownian motion, because the process $(\wh
W_s)_{s\geq 0}$ is not Markovian. We need to verify that, for every
choice of the rationals $u,v,u',v'$ such that $0<u<v<u'<v'$ we have
\begin{equation}
\label{techmini}
\min_{r\in[u,v]} \wh W_r\ne \min_{r\in[u',v']} \wh W_r,
\end{equation}
$\NN_0$ a.e. on the event $\{\sigma>v'\}$. In order to prove
(\ref{techmini}), we apply the Markov property at time $u'$ under
$\NN_0$. Notice that the law of $W_{u'}$ under $\NN_0(\cdot\mid
\sigma>u')$ is that of a Brownian path started from $0$ and stopped at
an independent random time whose law is explicitly known but
unimportant for what follows. So let us fix $\ell>0$ and denote by
$Q_\ell(d\w)$ the probability measure on $\W_0$ which is the law of a
linear Brownian motion started from $0$ and stopped at time $\ell$.
The proof of (\ref{techmini}) reduces to verifying that
\begin{equation}
\label{Markovmini}
\PP_\w\Big(\min_{r\in[0,v'-u']} \wh W_r=a,\; T>v'-u'\Big)=0\;,\
\forall a\in\RR\;,\ Q_\ell(d\w)\hbox{ a.s.}
\end{equation}
To simplify notation, set $h=v'-u'$. From the snake property, we have
for any $\w\in\W_0$
\begin{equation}
\label{mini0}
\min_{0\leq r\leq h}\wh W_r \leq
\Big(\min_{m(0,h)\leq t\leq \ell} \w(t)\Big)
\wedge \Big(\min_{m(0,h)\leq t\leq \zeta_h} W_h(t)\Big)
\end{equation}
$\PP_\w$ a.s.~on $\{T>h\}$. Hence, our claim (\ref{Markovmini}) will
follow if we can prove that $Q_\ell(d\w)$ a.s.,
\begin{equation}
\label{mini1}
\PP_\w\Big(\min_{0\leq r\leq h}\wh W_r=a,\;\min_{m(0,h)\leq t\leq
\ell} \w(t)>a,\;T>h\Big)=0
\;,\ \forall a\in\RR\;,
\end{equation}
and
\begin{equation}
\label{mini2}
\PP_\w\Big(\min_{0\leq r\leq h}\wh W_r=a=\min_{m(0,h)\leq t\leq \ell}
\w(t),\;T>h\Big)=0
\;,\ \forall a\in\RR\;.
\end{equation}

We first prove (\ref{mini1}), which in fact holds for every choice of
$\w$ with $\zeta_{(\w)}=\ell$, and not only $Q_\ell(d\w)$ a.s. We fix
$a\in \RR$. By properties of the Brownian snake, we know that under
$\PP_\w(\cdot\mid T>h)$ and conditionally on the pair
$(m(0,h),\zeta_h)$, the random path $(W_h(m(0,h)+t),0\leq t\leq
\zeta_h-m(0,h))$ is distributed as a linear Brownian motion started
from $\w(m(0,h))$ and stopped at time $\zeta_h-m(0,h)$. In particular,
on the event
$$
\Big\{\min_{0\leq r\leq h}\wh W_r=a,\;\min_{m(0,h)\leq t\leq \ell}
\w(t)>a,\;T>h\Big\}
$$
we have $\w(m(0,h))>a$ and thus
\begin{equation}
\label{mini3}
\min_{m(0,h)\leq t\leq \zeta_h} W_h(t)>a\ ,\ \PP_\w\hbox{ a.s.}
\end{equation}
because the law of the minimum of a Brownian path over a nontrivial
interval has a density, and we already know from (\ref{mini0}) that
the minimum in $(\ref{mini3})$ is greater than or equal to $a$ on the
event in consideration. We then argue by contradiction, assuming that
$$
\PP_\w\Big(\min_{0\leq r\leq h}\wh W_r=a,\;\min_{m(0,h)\leq t\leq
\ell} \w(t)>a,\;T>h\Big)>0\;.
$$
Using the Markov property under $\PP_\w$ at time $h$, together with
Lemma \ref{subtrees}, the property (\ref{mini3}) and the fact that
$\NN_x(\min_{s\geq 0}\wh W_s<y)<\infty$ for every $y<x$, it follows
that
$$
\PP_\w\Big(\min_{0\leq r\leq S_h}\wh W_r=a,\;\min_{m(0,h)\leq t\leq
\ell} \w(t)>a,\;T>h\Big)>0\;,
$$
where
$$
S_h:=\inf\{s\geq h:\zeta_s=m(0,h)\}\;.
$$
 From the definition of the ``excursions'' $W^i$ before Lemma
\ref{subtrees}, we then see that with positive probability under
$\PP_\w$, there exists some $i\in I$ such that
$$
\min_{s\geq 0} \wh W^i_s =a\ ,\ \w(\zeta_{\alpha_i})>a\;.
$$
However, by Lemma \ref{subtrees} and properties of Poisson measures,
the probability of the latter event is
$$
1-\exp\Big(-2\int_0^\ell dt\,{\bf 1}_{\{\w(t)>a\}}\,\NN_{\w(t)}
\Big(\min_{s\geq 0}\wh W_s=a\Big)\Big)=0
$$
because the law of $\min_{s\geq 0}\wh W_s$ under $\NN_x$ has no atoms
by (\ref{minisnake}). This contradiction completes the proof of
(\ref{mini1}).

It remains to prove (\ref{mini2}). We again fix $a\in \RR$. In
contrast with the previous argument, it will be important to disregard
certain sets of values of $\w$ which have zero $Q_\ell$-measure. We
first note that $Q_\ell(d\w)$ a.s.,
$$
\PP_\w\Big(\w(m(0,h))=\min_{m(0,h)\leq t\leq \ell}\w(t)\Big)=0
$$
so that the minimum of $\w$ over $[m(0,h),\ell]$ is attained $\PP_\w$
a.s.~at a point of $]m(0,h),\ell]$. It follows that $\PP_\w$ a.s.~on
the event
$$
\Big\{\min_{0\leq r\leq h}\wh W_r=a=\min_{m(0,h)\leq t\leq \ell}
\w(t),\;T>h\Big\}
$$
we can find a rational $q\in\;]0,\ell[$ such that, if $T_q=\inf\{s\geq
0:\zeta_s=q\}$,
$$
\min_{0\leq r\leq T_q} \wh W_r=a=\min_{q\leq t\leq \ell} \w(t)\;.
$$
Thus we need only check that the latter event has probability zero for
every rational $q\in\;]0,\ell[$. Using Lemma \ref{subtrees} once again,
we have
\begin{eqnarray*}
\PP_\w\Big(\min_{0\leq r\leq T_q} \wh W_r=
\min_{q\leq t\leq \ell} \w(t)\Big)
&=&\exp\Big(-2\int_q^\ell dr\,\NN_{\w(r)}\Big(\min_{s\geq 0} \wh
W_s\geq\min_{q\leq t\leq \ell} \w(t)\Big)\Big)\\
&=&\exp\Big(-3 \int_q^\ell dr\Big(\w(r)-\min_{q\leq t\leq \ell}
\w(t)\Big)^{-2}\Big)
\end{eqnarray*}
by (\ref{minisnake}). However, an application of L\'evy's modulus of
continuity for Brownian motion shows that
$$
\int_q^\ell dr\Big(\w(r)-\min_{q\leq t\leq \ell}
\w(t)\Big)^{-2}=\infty\ ,\ Q_\ell(d\w)\hbox{ a.s.}
$$
This completes the proof of (\ref{mini2}) and of Lemma
\ref{localminima}.
\cq

\smallskip
\noi{\bf Proof of Lemma \ref{keylemma}.}  %
We start by recalling that the law of the pair $(\eg_s,Z_s)_{0\leq
   s\leq 1}$ is invariant under time-reversal. More precisely,
$(\eg_{1-s},Z_{1-s})_{0\leq s\leq 1}$ has the same distribution as
$(\eg_s,Z_s)_{0\leq s\leq 1}$ (see e.g.~Section 2.4 in \cite{L}). Then
the proof of Lemma \ref{keylemma} reduces to checking that, almost
surely for every $s\in\;]0,1[$ such that
\begin{equation}
\label{key1}
\eg_s=\min_{r\in[s-\varepsilon,s]} \eg_r\quad,\ \hbox{ for some
}\varepsilon\in\;]0,s[\;,
\end{equation}
we have
\begin{equation}
\label{key2}
Z_s>\min_{r\in[s-\delta,s]} Z_r\quad,\
\hbox{ for every }\delta\in\;]0,s[
\end{equation}
and
\begin{equation}
\label{key3}
Z_s>\min_{r\in[s,s+\delta]} Z_r\quad,\
\hbox{ for every }\delta\in\;]0,1-s[\;.
\end{equation}
The fact that (\ref{key1}) implies (\ref{key2}) is an immediate
consequence of Lemma 2.2 in \cite{L} together with invariance under
time-reversal. We thus concentrate on the proof of (\ref{key3}). The
argument is similar to the proof of Lemma 2.2 in \cite{L}. We rely on
some ideas of Abraham and Werner \cite{AW}, which were already
exploited in Section 4 of \cite{LGW}.

Once again, we can reformulate the desired result in terms of the
Brownian snake.  It is enough to verify that, $\NN_0$ a.e.~for every
$s\in\;]0,\sigma[$ such that
$$
\zeta_s=\min_{r\in[s-\varepsilon,s]} \zeta_r\quad,\
\hbox{ for some }\varepsilon\in\;]0,s[\;,
$$
we have
$$
\wh W_s>\min_{r\in[s,s+\delta]} \wh W_r\quad,\ \hbox{ for every
}\delta\in\;]0,\sigma-s[\;.
$$

We first get rid of the case when $s$ corresponds to a local minimum
of the lifetime process $(\zeta_s)_{s\geq 0}$. Let $u$ and $v$ be two
rational numbers such that $0<u<v$, and argue under the probability
measure $\NN_0(\cdot\mid \sigma>v)$. We know that with probability one
there exists a unique time $s\in\,]u,v[$ such that $\zeta_s=m(u,v)$.
Furthermore, $\wh W_s=W_v(m(u,v))$ (by the snake property) and
conditionally on the pair $(m(u,v),\zeta_v)$, $(W_v(m(u,v)+t) -
W_v(m(u,v)),0\leq t\leq \zeta_v-m(u,v))$ is a linear Brownian path
started from $0$. From the snake property once again, we know that,
for every $\delta>0$, the set $\{\wh W_r, r\in [s,s+\delta]\}$
contains $\{W_v(t),m(u,v)\leq t\leq m(u,v)+\eta\}$ for some random
$\eta>0$ depending on $\delta$. By the preceding observations and
standard properties of linear Brownian paths, it follows that $\{\wh
W_r ,r\in[s,s+\delta]\}$ contains values strictly less than
$W_v(m(u,v))=\wh W_s$, as desired.

We now turn to the case when $s$ is not a time of local minimum of the
lifetime process.  Let us fix $u>0$, $\delta>0$ and an integer $A\geq
1$. It is enough to prove that, $\NN_0$ a.e. on the event
$\{u<\sigma\} \cap\{\zeta_u\leq A\}$, there exists no time
$s\in\;]u,\sigma[$ such that
\begin{equation}
\label{key4}
\zeta_s=\min_{r\in[u,s]} \zeta_r >2\delta
\end{equation}
and
\begin{equation}
\label{key5}
\wh W_r\geq \wh W_s\quad,\ \hbox{ for every }r\hbox{ such that }s\leq
r\leq \inf\{t>s:\zeta_t=\zeta_s-2\delta\}\;.
\end{equation}
To simplify notation, denote by $\NN^u_0$ the conditional probability
measure $\NN_0(\cdot\mid \sigma>u,\zeta_u\leq A)$. For every integer
$n\geq 1$ and every $i\in\{0,1,\ldots,A2^n\}$, set
$$
T^n_i=\inf\{r\geq u: \zeta_r=\zeta_u-i2^{-n}\}
$$
and
$$
S^n_i=\inf\{r>T^n_i:\zeta_r=\zeta_{T^n_i}-\delta\}\;,
$$
with the usual convention $\inf\emptyset=\infty$. Let $\alpha\in\;
]0,1]$. We can use the strong Markov property at time $T^n_i$,
together with Lemma \ref{subtrees} and (\ref{minisnake}), to evaluate
the probability
\begin{eqnarray}
\label{key6}
&&\NN^u_0(T^n_i<S^n_i<\infty;\wh W_r\geq \wh W_{T^n_i}-\alpha,\forall
r\in[T^n_i,S^n_i])=\\
&&\quad \NN^u_0\Big({\bf 1}{\{T^n_i<\infty,\zeta_{T^n_i}\geq \delta\}}\;
{\bf 1}{\{W_{T^n_i}(t) \geq \wh W_{T^n_i}-\alpha,\forall
t\in[\zeta_{T^n_i}-\delta,\zeta_{T^n_i}]\}}\nonumber\\
&&\qquad\qquad
\times\exp\Big(-2\int_{\zeta_{T^n_i}-\delta}^{\zeta_{T^n_i}} dt
\,\frac{3}{2} (W_{T^n_i}(t)-\wh W_{T^n_i}+\alpha)^{-2}\Big)\Big)\;.
\nonumber
\end{eqnarray}

By the snake property, on the event $\{T^n_i<\infty\}=\{\zeta_u\geq
i2^{-n}\}$, $W_{T^n_i}$ is the restriction of the path $W_u$ to the
time interval $[0,\zeta_u-i2^{-n}]$. Therefore under the probability
measure $\NN^u_0$ and conditionally on $\zeta_u$, the path $W_{T^n_i}$
is distributed as a linear Brownian path started from $0$ with
lifetime $\zeta_u-i2^{-n}$. It now follows that the quantity
(\ref{key6}) is equal to
\begin{eqnarray*}
&&\NN^u_0(T^n_i<\infty,\,\zeta_{T^n_i}\geq \delta)\;E\Big[{\bf
1}_{\{B_t\geq -\alpha,\forall t\in[0,\delta]\}}
\;\exp\Big(-3\int_0^\delta dt\,(B_t+\alpha)^{-2}\Big)\Big]\\
&&\quad\leq E\Big[{\bf 1}_{\{B_t\geq -\alpha,\forall t\in[0,\delta]\}}
\;\exp\Big(-3\int_0^\delta dt\,(B_t+\alpha)^{-2}\Big)\Big]
\end{eqnarray*}
where $(B_t,t\geq 0)$ is a linear Brownian motion starting from $0$
under the probability measure $P$. Finally we can use Proposition 2.6
in \cite{LGW} to get that the last quantity is bounded above by
$C_\delta\alpha^3$, where the constant $C_\delta$ only depends on
$\delta$ (compare with the estimate of Lemma 5.2 in \cite{L}).
Therefore we have obtained the bound
$$
\NN^u_0(T^n_i<S^n_i<\infty;\wh W_r\geq \wh W_{T^n_i}-\alpha,
\forall r\in[T^n_i,S^n_i])\leq C_\delta\,\alpha^3\;.
$$

We apply this estimate with $\alpha=2^{-2n/5}$. By summing over
possible values of $i$, and using the Borel-Cantelli lemma, we get
that, $\NN^u_0$ a.e., for every $n$ sufficiently large, for every
$i\in\{0,1,\ldots,A2^{n}\}$ such that $T^n_i<S^n_i<\infty$, the
condition
\begin{equation}
\label{key7}
\wh W_r\geq \wh W_{T^n_i}-2^{-2n/5}\hbox{ for every }r\in[T^n_i,S^n_i]
\end{equation}
does not hold.

To complete the proof, we argue by contradiction. Suppose that there
exists $s\in\;]u,\sigma[$ such that both (\ref{key4}) and (\ref{key5})
hold, and moreover $\zeta_u\leq A$. For every integer $n$ such that
$2^{-n}\leq \delta$, choose $i\in\{1,\ldots,A2^{n}\}$ such that
$T^n_{i-1}\leq s< T^n_i$. By (\ref{key4}), we have $\zeta_{T^n_{i-1}}
\geq \zeta_s$ and so $\zeta_{T^n_i}= \zeta_{T^n_{i-1}}-2^{-n}\geq
\zeta_s-\delta$. Hence,
$$
S^n_i\leq \inf\{r>s:\zeta_r=\zeta_s-2\delta\}\;,
$$
and by (\ref{key5}) we see that $\wh W_r\geq \wh W_s$ for every
$r\in[T^n_i,S^n_i]$. On the other hand, the snake property and
(\ref{key4}) ensure that $\wh W_s=W_u(\zeta_s)$, and we already
noticed that $\wh W_{T^n_i}=W_u(\zeta_{T^n_i})$. Since $0\leq
\zeta_s-\zeta_{T^n_i}\leq 2^{-n}$, the classical H\"older continuity
properties of Brownian paths imply that $\wh W_s\geq \wh
W_{T^n_i}-2^{-2n/5}$, for every $n$ sufficiently large. So we see that
for every $n$ sufficiently large, for $i$ chosen so that
$T^n_{i-1}\leq s< T^n_i$, the condition (\ref{key7}) holds. This
contradiction completes the proof of Lemma \ref{keylemma}.
\cq

\noindent \begin{tabular}{l}
D\'epartement de Math\'ematiques et Applications, UMR 8553 CNRS\\
Ecole Normale Sup\'erieure, 45 rue d'Ulm\\
75230 PARIS Cedex 05, FRANCE\\
{\it legall@dma.ens.fr, Frederic.Paulin@ens.fr}
\end{tabular}

\end{document}